\documentclass[a4paper, oneside,11pt]{amsart}
\usepackage{amssymb,mathrsfs,amstext, enumerate, comment}
\usepackage[noadjust]{cite}
\usepackage[plainpages=false,colorlinks,hyperindex,pdfpagemode=None,bookmarksopen,linkcolor=blue,citecolor=blue,urlcolor=blue]{hyperref}

 \usepackage{bbm}
 \usepackage{graphicx}

\newtheorem{thm}{Theorem}[section]
\newtheorem{lem}[thm]{Lemma}

\newtheorem{prop}[thm]{Proposition}
\newtheorem{cor}[thm]{Corollary}
\newtheorem{defn}[thm]{Definition}

\numberwithin{equation}{section}

\def\H{{\mathscr{H} }}

\def\M{{\mathscr{ M}}}

\def\M{{\mathscr{M} }}

\def\a{{\alpha}}

\begin{document}

\title[$W^{\ast}$-TROs]{The linking von Neumann algebras of $W^{\ast}$-TROs}

\author[Wang]{  Liguang Wang}
\address[Wang]{School of Mathematical Sciences, Qufu normal University, Qufu 273165, China.}
\email{wangliguang0510@163.com}

\author[Chen]{Hongjie Chen}
\address[Chen]{School of Mathematical Sciences, Qufu Normal University, Qufu, Shandong, 273165, China}
\email{chenhongjiehei@163.com}

\author[Wong]{Ngai-Ching Wong}
\address[Wong]{Department of Applied Mathematics, National Sun Yat-sen University, Kaohsiung, 80424, Taiwan, \and
	School of Mathematical Sciences, Tiangong University, Tianjin 300387, China.}
\email{wong@math.nsysu.edu.tw}

% \thanks{Liguang Wang is the corresponding author.}

\thanks{Wang is partially supported in part by NSF of China (NSF of China (12471124)).}

\thanks{Wong is grateful to the colleagues in Tiangong University at which he spent a fruitful sabbatical leave recently.}

\thanks{Corresponding author: Ngai-Ching Wong, wong@math.nsysu.edu.tw}

\begin{abstract}
In this  note, we  show that a von Neumann algebra
can be written as the linking von Neumann algebra of a $W^\ast$-TRO
if and only if it contains no  abelian direct summand.
We also provide some new characterizations of
 nuclear TROs and $W^\ast$-exact TROs in terms of the properties of their linking algebras.
  \end{abstract}

 \maketitle

 \parskip=\baselineskip
 %\vskip0.1cm

\section{Introduction}

A \emph{ternary ring of operators} (or simply \emph{TRO}) %(\cite{H})
is a norm closed subspace $V$ of the Banach space $B(K, H)$ of bounded linear operators between  Hilbert spaces $K$ and $H$,
which is closed under the triple product
$$
(x, y, z)\in V\times V\times V\mapsto xy^{\ast}z\in V.
$$
TROS were first introduced by Hestenes (\cite{H}), and pursued by many others.
% (see e.g. \cite{Har,Z,Ha1, Ha2,EOR, Kaur,Ruan2}).
In \cite{NR},
it is proved that TROs form a special class of concrete operator spaces and
characterized TROs in terms of the operator space theoretic properties. The interconnections between TROs and JC$^{\ast}$-triples are studied in \cite{BT}.

When $V$ is a TRO, $V^{\sharp}=\{x^{\ast}\in B(H, K): x\in V\}$ is also a TRO.
We   assume that $V$ is \emph{non-degenerate} in this note, in the sense
that $VK$ and $V^\sharp H$ are norm dense in $H$ and $K$, respectively.
A TRO $V$ of $B(K, H)$ is called a $W^\ast$-TRO if it is closed in the strongly operator topology (SOT), or equivalently,
closed in the weak operator topology, or the  weak$^{\ast}$ topology of  $B(K, H)$ (\cite{Z}).
%\end{defn}

%Ternary rings of operators
% were first introduced by Hestenes \cite{H} and have been intensively studied by Harris \cite{Har}, Zettl \cite{Z},
%Hamana \cite{Ha1, Ha2},  %Effors-Ozawa-Ruan \cite{EOR}, Kaur and Ruan \cite{Kaur} and Ruan \cite{Ruan2}.

A fundamental tool to study TROs
is the construction of the linking algebra, that is, a particular C$^{\ast}$-algebra containing
the related TRO as a corner.  For example,
an operator space is injective if and only if it is completely isometric to a ternary
corner of an injective C$^{\ast}$-algebra (see, e.g., \cite{BL}).
Let $VV^{\sharp}$ and $V^{\sharp}V$ be the linear span of $vw^{\ast}$ and $v^{\ast}w$ for all $v, w\in V$ respectively.
Clearly, $VV^{\sharp}$ and $V^{\sharp}V$  are $\ast$-subalgebras of $B(H) $ and $B(K)$.
Let
$$
C(V)=\overline{VV^{\sharp}}^{||\cdot||}\quad\text{and}\quad D(V)=\overline{V^{\sharp}V}^{||\cdot||}
$$
denote the $C^\ast$-algebras generated by $VV^{\sharp}$ and $V^{\sharp}V$ respectively.
The \emph{linking $C^*$-algebra} $A(V)$ of $V$ is defined by
$$
A(V)=\left(
                                              \begin{array}{cc}
                                                C(V) & V \\
                                                 V^{\sharp}   & D(V)\\
                                              \end{array}
                                            \right).
$$
When $V$ is a $W^\ast$-TRO,  let
$$
M(V)=\overline{V V^{\sharp}}^{\operatorname{SOT}}\quad\text{and}\quad N(V)=\overline{V^{\sharp}V}^{\operatorname{SOT}}
$$
 denote  the von Neumann algebras generated by $VV^{\sharp} $ and $V^{\sharp}V$ respectively.
The \emph{linking von Neumann algebra} $R(V)$ of $V$ is defined by
$$
R(V)=\left(
                                              \begin{array}{cc}
                                                M(V) & V \\
                                                 V^{\sharp}   & N(V)\\
                                              \end{array}
                                            \right)
=A(V)'',
$$
the double commutant of $A(V)$ in $B(H\oplus K)$ (see, e.g., \cite[Proposition 2.3]{EOR}).

TROs and their associated linking algebras share many
common properties, wherefore the application of operator algebraic methods simplifies
the study of TROs that are not algebras themselves. Basic properties and most recent results of TROs are
discussed in, e.g., \cite{EM, GJL, H, Har, KK, KKR, Ha1, Ha2, EOR, Kaur, Ruan2, SV, Z} and references therein.

In this note,
%we  consider when a von Neumann algebra can be the linking von Neumann algebra of a $W^\ast$-TRO.
 we show that a von Neumann algebra  $\M$ can be written as  the linking von Neumann algebra $R(V)$
  of a $W^\ast$-TRO $V$ if and only if $\M$ does not contain an
 abelian direct summand.
We also provide new characterizations of a TRO being nuclear or $W^\ast$-exact in term of its linking algebra.
 These results  generalize, in particular,  \cite[Theorem 6.5]{Kaur}  and \cite[Theorem 4.1]{DR}.

\section{The results}

\subsection{Conditions to be a linking von Neumann algebra of a $W^*$-TRO}

\begin{thm} \label{thm:Ce=C(1-e)=1}
Let $\M$ be a von Neumann algebra.
The following conditions are equivalent.
\begin{enumerate}[(a)]
  \item There is a $W^\ast$-TRO $V$ such that its linking von Neumann algebra $R(V)=\M$.
  \item There exists a    projection $ e$ in $\M$ with central covers   $C_e= C_{I-e}=I$.
  \item $\M$ has no  abelian direct summand.
\end{enumerate}
\end{thm}
\begin{proof}
For the implication (a) $\implies$ (b), suppose
 $\M=R(V)$ is the linking algebra of  $W^\ast$-TRO $V$  as in \eqref{eq:RV}.
Then $e=\begin{pmatrix} 1&0\\0&0\end{pmatrix}$ is a projection in $\M$ with $C_e = C_{I-e}=I$.

For the implication (b) $\implies$ (a), suppose there exists a    projection $ e$ in $\M$
with central covers   $C_e= C_{I-e}=I$.
Let $V=e\M (I-e)$. Then $V$ is a $W^\ast$-TRO with
\begin{align}\label{eq:RV}
R(V)=\left(
                                              \begin{array}{cc}
                                                M(V) & V \\
                                                 V^{\sharp}   & N(V)\\
                                              \end{array}
                                            \right).
\end{align}
 %where $M(V)=\overline{VV^{\sharp}}^{s.o.t}$ and $N(V)=\overline{V^{\sharp }V}^{s.o.t}$.
We claim that $R(V)=\M$.

To see $M(V)=e \M e$, it suffices to show that $e \M_+ e\subset M(V)$.
 Since
 $$
I= C_{I-e}=\bigvee \{ {u}(I-e){u}^{\ast}:  u\ \text{ is a unitary in $\M$}\},
 $$
 there is a net $\{F_i\}_{i}$ of increasing finite
 subsets of the unitary group  of $\M$ such that
 $$
 \bigvee_{{u}\in F_i}{u} (I-e){u}^{\ast}\rightarrow I
 $$
  in the
 strong operator topology.
 Given $x\in \M_+$.
 Since
$ex^{1/2}{u} (I-e) {u}^{\ast}x^{1/2}e$ belongs to the von Neumann algebra  $M(V)$, we have
$$
exe= \operatorname{SOT-}\lim_i ex^{1/2}\left(\bigvee_{{u}\in F_i}{u} (I-e){u}^{\ast}\right)x^{1/2}e\in M(V).
$$

Similarly, we see that $N(V) = (I-e)\M (1-e)$.  It follows $\M=R(V)$.

For the implication (c) $\implies$ (b),
suppose $\M$ has no
 abelian direct summand.
We show that
there is a    projection $e$ in $\M$ such that the central covers   $C_e=C_{I-e}=I$.
Indeed, the assertion is contained in \cite[Exercise 6.19]{KR}.   For completeness, we present a short proof below.

Write $\M =\M_d \oplus \M_c$ as the direct sum of its discrete part $\M_d$ and continuous part $\M_c$, with identity elements
$I_{\M_d}$ and $I_{\M_c}$, respectively.
For a discrete summand
 $M_{n} (\mathbb{C})\otimes L_{\infty}(\mu_n)$ (here $n\geq 2$), let $e_{n}=E^n_{11}\otimes I_{L_{\infty}(\mu_n)}$,
 where $E^n_{11}$ is the matrix unit in $M_{n}(\mathbb{C})$ with the $(1,1)$ entry being $1$ and all others $0$.
Let $e_d = \sum_{n\geq 2} e_{n}$.   Then $C_{e_d} = C_{I_{\M_d} - e_d} = I_{\M_d}$.
For the continuous part, we have a projection $e_c \in \M_c$ such that $e_c$ and $I_{\M_c}-e_c$ are both unitarily equivalent to $I_{\M_c}$.
In particular, $C_{e_c}=C_{I_{\M_c} - e_c} = I_{\M_c}$.  Consequently, $e = e_d + e_c$ finishes our task.

Finally, we verify the implication (b) $\implies$ (c).
Let $e$ be  a projection in $\M$ such that $C_e=C_{I-e}=I$.
Let $z$ be any abelian central projection in $\M$.
 Since
 $$
 ze(u^*(1-e)u)ez = zez(zuz)^*(1-e)uez=(zuz)^*(zez)(1-e)u=0
 $$
for any unitary $u$ in $\M$, and $C_{I-e}=I$, we see that $ze=0$.  Similarly, we see that $z(1-e)=0$, and thus $z=0$.
It follows that $\M$ has no direct abelian summand.
\end{proof}

%It is clear that the commutative von Neumann algebra $\mathbb{C}$ cannot be the linking algebra of any $W^\ast$-TRO.

\subsection{Characterization of nuclearity}
Nuclear $C^\ast$-algebras play  an important role in the study of operator algebras.
Analogously, nuclear TROs
are also characterized in \cite{Kaur}.
Our aim   is to give some new characterizations of nuclear TROs.

Recall that   a $C^\ast$-algebra $A$ (resp.\ TRO $V$) is said to be \emph{nuclear} (resp.\ \emph{Lance nuclear}) if for every
$C^\ast$-algebra $B$, there is a unique $C^\ast$-algebra tensor norm on $A\otimes B$ (resp.\ a unique TRO tensor norm on
$V\otimes B$).

\begin{prop}\label{3.1} If $V$ is a TRO in $B(K,H)$ such that $C(V)$ and $D(V)$ are   nuclear,
then $A(V)$ is nuclear and  $V$ is Lance nuclear.\end{prop}

\begin{proof}  There is a projection $e$ in $B(H\oplus K)$ such that $C(V)=eA(V)e$ and $D(V)=(I-e)A(V)(I-e)$.
Hence $C(V)^{\ast\ast}=eA(V)^{\ast\ast}e$ and $D(V)^{\ast\ast}=(I-e)A(V)^{\ast\ast}(I-e)$. Since $C(V)$ and $D(V)$ are  nuclear,
 $C(V)^{\ast\ast}$ and $D(V)^{\ast\ast}$ are hyperfinite. We have  $A(V)^{\ast\ast}$ is hyperfinite by
   \cite[Lemma 2.8]{Wang2}. Hence $A(V)$ is   nuclear, and therefore $V$ is Lance nuclear by  \cite[Theorem 6.1]{Kaur}.
\end{proof}

A  TRO (respectively $W^\ast$-TRO) $V\subset B(K, H)$ carries
 a natural operator space structure (\cite{ER}; see also  \cite{P}) with
 matrix norms arising from identifying $M_n(V)$ with a TRO (respectively $W^\ast$-TRO)   in $M_n(B(K, H))=B(K^n, H^n)$ for   $n=1,2,\ldots$.
In general, an operator space $X$ is said to be \emph{injective}
if for any operator spaces $W_1 \subseteq W_2$, every complete contraction
$\phi : W_1 \to X$ has a completely contractive extension $\hat{\phi} : W_2 \to X$. % with $\|\hat{\phi}\|_{\mathrm{cb}}\leq 1$.
On the other hand, $X$ is said to be
\emph{$1$-nuclear} if the identity operator $I_X$ on   $X$
can be factorized approximately through   matrix spaces $M_{n_\alpha}(\mathbb{C})$
in the sense that $I_X = \lim_\alpha \phi_\alpha\circ \psi_\alpha$ in the  point-norm topology with completely bounded maps
$\psi_\alpha: X\to M_{n_\alpha}(\mathbb{C})$ and $\phi_\alpha : M_{n_\alpha}(\mathbb{C})\to X$
such that $\|\phi_\alpha\|_{\mathrm{cb}}\|\psi_\alpha\|_{\mathrm{cb}}\leq 1$.

Arguing as in Proposition \ref{3.1}, we have a similar result held  for $W^\ast$-TROs.

\begin{prop}  If $V$ is a $W^\ast$-TRO such that $M(V)$ and $N(V)$ are injective von Neumann algebras, then $R(V)$ is an injective von Neumann algebra and therefore  $V$ is an injective TRO.\end{prop}

It follows from Proposition \ref{3.1}  the following characterization of nuclear TROs, which
adds condition (6) to the list in \cite[Theorem 6.5]{Kaur}

%  We  refer to  \cite{Kaur} for related notations.

\begin{thm}\label{Thm3.2}
Let $V$ be a TRO. The following are equivalent:
\begin{enumerate}
\item[(1)] $V$ is Lance nuclear;

\item[(2)] $V$ is 1-nuclear;

\item[(3)] $V^{\ast\ast}$ is injective;

\item[(4)] $A(V)^{\ast\ast}$ is injective;

\item[(5)] $A(V)$ is nuclear;

\item[(6)] $C(V)$ and $D(V)$ are nuclear.
\end{enumerate}
\end{thm}

\subsection{Characterization of exactness}
Recall that a von Neumann algebra $\M$ is said to be \emph{weakly exact} (\cite{Ozawa}) if for any
unital $C^\ast$-algebra $A$ with a
 closed two-sided ideal $J$  and
any left normal $*$-representation $\pi: \M\otimes_{\min} A\to B(H)$ with $\pi(\M\otimes J)=0$, the induced $*$-representation
$\hat{\pi}: \M\otimes (A/J)\to B(H)$ is continuous with respect to the minimum tensor norm.

%The following  result is   Corollary 14.1.15 in Brown-Ozawa's book (\cite{BO}).

\begin{lem}[{\cite[Corollary 14.1.15]{BO}}] \label{AAA}
If $\M$ is a weakly exact von Neumann algebra, then $\M\overline{\otimes} B(\H)$
 is weakly exact.%, in particular, $ \M\overline{\otimes }M_2(\mathbb{C})$ is weakly exact.
\end{lem}

Dong and Ruan studied the connection between weak$^{\ast}$ exact $W^\ast$-TROs and their linking von Neumann algebras in \cite{DR}.
In view of \cite[Theorem 3.3]{DR}, a dual operator space $X$ is \emph{weak* exact} if
for any operator space $W$ and any finite rank complete contraction $\phi : W\to X$,
 there exists
a net of weak* continuous finite rank complete contractions $\phi_\alpha: W\to V$   converging
to $\phi$ in the point-weak* topology.

\begin{prop}[{\cite[Lemma 3.4]{DR}}] \label{DR3.6}
 Let $\M$ be a von Neumann algebra. Then $\M$ is weak$^{\ast}$ exact if and only if $\M$ is weakly exact.
 \end{prop}

\begin{comment}
\begin{defn}[{\cite{DR}}] \label{defn:C}
A TRO $V $ is said to satisfy \emph{condition $C^{'}$}
 if for any TRO $W$, there is a completely isometric isomorphism
$$ V\overset{\vee }{\otimes}: W^{\ast\ast}=V\overset{\vee }{\otimes}   W^{\ast\ast}.$$
A dual TRO $V$ is said to satisfy \emph{weak$^{\ast}$ condition $C^{'}$} or
to be \emph{weak$^{\ast}$-exact} if for   any TRO $W$, there is a completely isometric isomorphism
$$ V\overset{\sigma,\vee }{\otimes}:   W^{\ast\ast}=V\overset{ \vee }{\otimes}   W^{\ast\ast}$$
\end{defn}
\end{comment}
%\textcolor[rgb]{1.00,0.00,0.00}{Wong:  I do not understand this definition, neither the original or this new form.
%Please explain clearly, and give a reference if it is copied from others.  Moreover, since the notion of operator space is not
%defined, and indeed not discussed anywhere, it would be better to make the definition works for TROs.}

%Q{The following results  are needed later.}

\begin{lem}[{\cite[ Lemma 3.1]{DR}}]\label{DR3.5}
Assume that $X_0$ is a weak$^{\ast}$-closed subspace of a dual operator space $X$
such that there exists a weak$^{\ast}$ continuous completely contractive  projection $P: X\rightarrow X_0$. If $X$
is weak$^{\ast}$ exact, so is $X_0$.
\end{lem}

{The following result connects the weak exactness of a von Neumann algebras with its ``diagonals''.}

\begin{lem}\label{lem:diagonal}
Let $\M$ be a  von Neumann algebra. Suppose $\{Q_{j}: j \in \mathbb{J}\}$ is a family of
mutually equivalent and mutually orthogonal projections in $\M$ such that
$Q_j\M Q_j$ is a weakly exact von Neumann algebra for each $j$
and $\sum_{j\in \mathbb{J}} Q_j=I$. Then $\M$ is a weakly exact von Neumann algebra.
\end{lem}

\begin{proof}
Fix an index $j_0$ in $\mathbb{J}$.
 For any finite subset $F$ of $\mathbb{J}$ of $n$ elements,
 we see that
 \begin{quote}
  $(\sum_{j\in F}Q_j)\M (\sum_{j\in F}Q_j)$   is isomorphic to    $\left(Q_{j_0}\M Q_{j_0}\right)\otimes M_n(\mathbb{C})$.
  \end{quote}
 But the latter von Neumann algebra is weakly exact by Lemma \ref{AAA}. Let $P_F(x)=(\sum_{j\in F}Q_j)x(\sum_{j\in F} Q_j)$.
 Then as $F$ runs over all finite subsets of $\mathbb{J}$, we have a net of normal completely positive
  contractions $P_{F}: \M\rightarrow P_F\M P_F$ which converges to the identity on $\M$
  in the point ultraweak topology. It follows from
  \cite[Proposition 14.1.4]{BO}  that $\M$ is a weakly exact von Neumann algebra.
\end{proof}

\begin{lem}[{\cite[Theorem 3.1]{DYH}; see also \cite[Lemma 2.7]{Wang2}}]\label{ABC}
Let $P$ be a projection in  a von Neumann algebra   $\M $.
There is a family $\{Q_{\alpha} \}_\alpha$ of subprojections of $P$ in $\M $ such that
$C_P = \sum_{\alpha}C_{Q_\alpha}$ is a sum of mutually orthogonal central projections. Moreover, each
$C_{Q_\alpha} = Q_\alpha + \sum_{i}Q_{\alpha}^{i}$
is a sum of mutually orthogonal and mutually equivalent projections.
\end{lem}
\begin{proof}
  We sketch the proof in \cite{DYH} for completeness.
  We first claim that for any nonzero projection $P$ in $\M$ there is a subprojection $Q\leq P$ and a family
%$\mathcal{Q}$
of mutually orthogonal projections $Q_\alpha\sim Q$ such that the central cover $C_Q = Q + \sum_\alpha Q_\alpha$.
The assertion will follow  from the claim and a Zorn's Lemma argument.

To prove the claim, we might assume $C_P=I$.
Suppose $\M$ has type III.
There are subprojections $P_1, P_2$ of $P$ such that $P=P_1 + P_2$ and $P_1\sim P_2$.
Enlarge $\{P_1, P_2\}$ to a maximal family $\mathcal{P}$   of mutually orthogonal
projections $P_\alpha$ such that $P_\alpha\sim P_1$ for all indices $\alpha$.
If $\sum_\alpha P_\alpha=I$, then we can let $Q=P_1$ and $\{Q_{\alpha} \}_\alpha=\mathcal{P}\setminus\{P_1\}$.
Otherwise, $P_1\not\precsim (1- \sum_{\alpha\neq 1} P_\alpha)$ by the maximality of $\mathcal{P}$.
Hence there is a nonzero central projection $E$ of $\M$ such that $E(1- \sum_{\lambda\neq 1} P_\lambda) \prec EP_1$.
Then $Q=EP_1$ and $Q_\alpha = EP_\alpha$ for $\alpha\neq 1$ will do the job.

If $\M$ is semifinite, there will be a nonzero central projection $E=\sum_{\alpha\in \Lambda} E_\alpha$ written
as an orthogonal sum of mutually equivalent finite projections $E_\alpha$.  Since $C_P=I$, we can replace
$P$ by $EP\neq 0$ and assume $I=E=\sum_{\alpha\in \Lambda} E_\alpha$.
Letting $E_{\alpha\beta}=E_{\beta\alpha}^*$
be the partial isometry in $\M$ such that $E_\alpha =E_{\alpha\beta}E_{\alpha\beta}^*$ and
$E_\beta =E_{\alpha\beta}^*E_{\alpha\beta}$,
we get a family $\{E_{\alpha\beta}: \alpha, \beta\in \Lambda\}$ of matrix units
 in $\M$ with $E_{\alpha\alpha}=E_\alpha$
for each $\alpha$ in $\Lambda$.  Then
$U_{\alpha\beta} = E_{\alpha\beta}+E_{\beta\alpha} +
\sum_{\gamma\neq \alpha,  \beta} E_{\gamma\gamma}=U_{\beta\alpha}$
is a self-adjoint
unitary in $\M$ such that $U_{\alpha\beta}E_{\alpha}=E_{\beta}$, $U_{\alpha\beta}E_{\beta} = E_{\alpha}$ and
$U_{\alpha\beta}E_{\gamma} =E_{\gamma}$ for $\gamma\neq \alpha, \beta$.

If there is a nonzero central projection $F$ in $\M$ such that $FE_\alpha \precsim FP$
for some (and thus all) $\alpha$,
then let $Q\leq FP$ such that $Q\sim FE_\alpha$.
Since $FE_\alpha$ is a finite projection, there is a unitary $U$ in $\M$ such that
$Q=U^*FE_\alpha U$.  Letting $Q_\beta = U^*FE_\beta U$ for all $\beta\neq \alpha$ in $\Lambda$, we have
$C_Q=F = Q + \sum_{\beta\neq \alpha} Q_\beta$ as claimed.

If there is no such nonzero central projection $F$, then $P\prec E_\alpha$ for all $\alpha$ in $\Lambda$.
In particular, $P\leq U^*E_{\alpha} U$ for some index $\alpha$ and a unitary $U$ in $\M$.
Replacing all $E_\beta$ with $U^*E_\beta U$, we can assume $U=I$,
and thus $P$ is  a projection in the finite von Neumann
 algebra $E_\alpha \M E_\alpha$.
   By \cite[Proposition 8.2.1]{KR}, $P$ contains a nonzero monic subprojection $Q$; namely,
 there are mutually orthogonal subprojections $Q_1=Q, Q_2 \ldots, Q_k$ of $E_\alpha$, each of them is equivalent to $Q$ such
 that $E_{\alpha} C_Q E_{\alpha}= Q_1 + \cdots + Q_k$.
 Consequently,
\begin{align*}
C_Q &=  \sum_{\beta} E_{\beta}C_Q E_{\beta}= E_{\alpha} C_Q E_{\alpha} +
\sum_{\beta\neq\alpha} E_{\beta}C_Q E_{\beta} \\
&= (Q_1 + \cdots + Q_k) + \sum_{\beta\neq\alpha} U_{\alpha\beta}(Q_1 + \cdots + Q_k) U_{\alpha\beta}.
\end{align*}
Because each $U_{\alpha\beta}Q_j U_{\alpha\beta}$ is equivalent to $Q$, we establish the claim.

Since every von Neumann algebra is a direct sum of its semifinite summand and its type III summand, the claim is verified.
\end{proof}

\begin{prop} \label{prop3.9}
Let $\M$ be a  von Neumann algebra and $e\in \M$ a nonzero projection with $C_{e}=I$.
Then
\begin{center}
$\M$ is weak$^{\ast}$-exact \quad  $\Longleftrightarrow$ \quad  $e\M e$  is
weak$^{\ast}$-exact.
\end{center}
\end{prop}

 \begin{proof}  ($\Rightarrow$) It follows from Lemma \ref{DR3.5} that if $\M$ is weak$^{\ast}$-exact, then $e\M e$ and $(I-e)\M(I-e)$ are weak$^{\ast}$-exact.  This is because
 the maps $P(x)=exe$ and $Q(x)=(I-e)x(I-e)$ are weak$^{\ast}$-continuous completely contractive
   projections from $\M$ onto $e\M e$ and $(I-e)\M (I-e)$, respectively.

($\Leftarrow$)   % We only consider the case when
Assume that $e\M e$  is
weak$^{\ast}$-exact.
By Lemma \ref{ABC},   there is a family $\{Q_{\alpha}: \alpha\in  {\Gamma}\}$ of mutually disjoint subprojections of  $e$
such that $\sum_\alpha C_{Q_{\alpha}}=1$, and
$C_{Q_{\a}}=Q_{\a}+\sum_i Q_{\a}^{i}$ where all $Q_{\a}^{i}$ are equivalent to $Q_{\a}$ for each $\alpha$.
By Lemmas \ref{DR3.6} and \ref{DR3.5}, we have $Q_{\a}\M Q_{\a}$ is weakly exact for each $\a\in {\Gamma}$.
 It follows from
 Lemma \ref{lem:diagonal} that $C_{Q_{\a}}\M C_{Q_{\a}}$ is a weakly exact von Neumann algebra.
Since the direct sum of weakly exact von Neumann algebras is again weakly exact,
 $\M=\sum_{\a\in {\Gamma}} C_{Q_{\a}}\M $ is a weakly exact von Neumann algebra. This completes the proof.
\end{proof}

%\textcolor[rgb]{1.00,0.00,0.00}{Wong:  I do not understand the  proof.  I guess it is claimed that
%if $C_e = C_{I-e}=I$ then there is a mutually orthogonal and equivalent family of projections $p_\alpha$
%summing to $I$ such that each $p_\alpha\M p_\alpha$ is weakly exact.  But I do not see how to get this family.
%For example $e=E_{11}$ in $M_3(\mathbb{C})$.  In the original argument, $C_e = I$ for any nonzero projection $e$.
%Maybe there is a typo, and one should red $\sum_\alpha Q_\alpha =I$ instead of $\sum_\alpha C_{Q_\alpha} =I$.}\\

%{\color{blue} I provide Lemma 3.9 after Lemma 3.8.}\\

It follows from Proposition \ref{prop3.9} that
for a   $W^\ast$-TRO $V$, the linking $W^*$-algebra $R(V)$ is weak$^{\ast}$-exact if and only if
$M(V)$  or $N(V)$ is a weakly exact von Neumann algebra.
The following result is immediate and it adds condition (2') to the list in  \cite[Theorem 4.1]{DR}.

\begin{thm}  Let $V$ be a $W^\ast$-TRO. Then the following are equivalent.
\begin{enumerate}
\item[(1)] $V$ is weak$^{\ast}$-exact.

\item[(2)] $M(V)$ and $N(V)$ are weak$^{\ast}$ exact.

\item[(2')] $M(V)$ or $N(V)$ is weak$^{\ast}$ exact.

\item[(3)] $R(V)$ is weak$^{\ast}$ exact.
\end{enumerate}
\end{thm}

\begin{comment}
It follows from Theorem \ref{lem:Linking-alg} the following result.

\begin{cor}
  Let $\M$ be a von Neumann algebra without abelian direct summand.  Then $\M$ is weakly exact if and only if $e\M e$ is
  weakly exact for any, and thus all, projection $e$ in $\M$ with central covers $C_e = C_{I-e} = I$.
\end{cor}
\end{comment}

\end{document}